# ON FINITE GROUPS WITH POLYNOMIAL DIAMETER

LUCA SABATINI

ABSTRACT. Given a finite group $G$ and a generating set $S \subseteq G$, the diameter $diam(G, S)$ is the least integer $n$ such that every element of $G$ is the product of at most $n$ elements of $S$. In this paper, for bounded $|S|$, we characterize groups with polynomial diameter as the groups with a large abelian section close to the top, precisely of size an exponential portion of the size of the full group. This complements a key result of Breuillard and Tointon [2]. As a consequence, groups with polynomial diameter have many conjugacy classes, and contain a large nilpotent subgroup of class at most 2.

## 1. INTRODUCTION

Let $G$ be a finite group, and $S \subseteq G$. The diameter $diam(G, S)$ is the least integer $n$ such that every element of $G$ is the product of at most $n$ elements of $S$ (if $S$ does not generate $G$, we assume that this value is infinite). It is quite intuitive that a large abelian portion of the group is related to a large diameter: in this short paper, we make this insight very clear.

**Theorem 1.1** (Structure of groups with polynomial diameter). *Let $G$ be a finite group, and $\varepsilon, c > 0$. The following hold:*
  (A) *If there exists a generating set $S \subseteq G$ such that $diam(G, S) \geq c|G|^\varepsilon$, then there exists $H \leqslant G$ such that*
  $$|G : H| \ll_{\varepsilon, c, |S|} 1 \qquad \text{and} \qquad |H/H'| \gg_{\varepsilon, c, |S|} |G|^\varepsilon.$$
  (B) *If there exists $H \leqslant G$ such that $|H/H'| \geq c|G|^\varepsilon$, then for every set $S \subseteq G$ one has*
  $$diam(G, S) \geq c' \cdot |G|^{\frac{\varepsilon}{|G:H||S|}},$$
  *where $c' := \min(1, c)$.*

The function $\varepsilon |G : H|^{-1}|S|^{-1}$ in (B) is the best possible. The "large abelian section" condition does not involve the generating set anyway, and it is purely related to the structure of $G$. Hence, the polynomial diameter is a *group property*, and does not depend on the choice of the generating set (this fact is also proved in [2, Corollary 4.15]). Moreover, we are able to provide an estimate for the diameter, which is interesting also for sets of unbounded size.

**Proposition 1.2** (Changing the set). *Let $diam(G, S) \geq c|G|^\varepsilon$ for some $\varepsilon, c > 0$ and generating set $S$. Let $S' \subseteq G$ be another set. Then there exists $\varepsilon' = \varepsilon'(\varepsilon, c, |S|) > 0$ such that*
$$diam(G, S') \gg_{\varepsilon, c, |S|} |G|^{\varepsilon'/|S'|}.$$

*Proof (from Theorem 1.1).* Applying Theorem 1.1 (A) and then Theorem 1.1 (B), we obtain the desired result with $\varepsilon' := \frac{\varepsilon}{|G:H|}$, where $|G : H| \ll_{\varepsilon, c, |S|} 1$. □

As byproducts, we prove other interesting properties of groups with large diameter.

**Proposition 1.3.** *If $diam(G, S) \geq c|G|^\varepsilon$ for some generating set $S \subseteq G$, then*
  (i) *$G$ has at least $O_{\varepsilon, c, |S|}(|G|^\varepsilon)$ conjugacy classes;*
  (ii) *$G$ has a nilpotent subgroup of size at least $O_{\varepsilon, c, |S|}(|G|^\varepsilon)$ and class at most $2$.*

## 2. PRELIMINARIES

By $f(n) \ll g(n)$, $g(n) \gg f(n)$ and $f(n) = O(g(n))$, we mean the same thing, namely, that there is $C > 0$ such that $|f(n)| \leq C \cdot g(n)$ for all $n \geq C$. The constant $C$ is allowed to depend on the variables in the subscript of $\ll, \gg, O$.

Let $H \leqslant G$. As it is well described in [3, Section 7.2], one can derive a set of generators for $H$ from a set of generators of $G$, by the usual Reidemeister-Schreier process.

**Lemma 2.1** (Lemma 7.2.2 in [3])**.** *Let $G$ be a finite group, $H \leqslant G$, and $S \subseteq G$. Then there exists $\overline{S} \subseteq H$ of size at most $|G:H||S|$ such that $diam(H, \overline{S}) \leq diam(G, S)$.*

The following result is folklore: we give a proof for the convenience of the reader.

**Lemma 2.2** (Diameter of abelian groups)**.** *Let $G$ be a finite **abelian** group, and $S \subseteq G$. Then*
$$diam(G, S) \geq |G|^{1/|S|}.$$

*Proof.* Let $n := diam(G, S)$, and let $s_1, ..., s_{|S|}$ be the elements of $S$. Since $G$ is abelian, the set of words of length at most $n$ in $S$ is precisely $(s_1)^{\alpha_1} \cdot ... \cdot (s_{|S|})^{\alpha_{|S|}}$, where $\sum_i \alpha_i \leq n$. Of course, these words are at most $n^{|S|}$. It follows that $|G| \leq n^{|S|}$, as desired. $\square$

For the proof of Theorem 1.1 (A), a crucial role is played by a result of Breuillard and Tointon concerning product sets. When $S \subseteq G$ and $n \in \mathbb{Z}_+$, we write
$$S^n := \{s_1 \cdot ... \cdot s_n : s_1, ..., s_n \in S\}.$$
The following is implicitly obtained along the proof of [2, Theorem 4.1 (1)-(2)].

**Theorem 2.3** (Breuillard-Tointon)**.** *Let $\theta, \delta > 0$, and $S \subseteq G$ a generating set such that*
$$|S^{5n}| \leq \theta \cdot |S^n|$$
*for some $n \leq diam(G, S)^\delta$. Then there exists a normal subgroup $N \triangleleft G$ contained in $S^{\lfloor diam(G,S)^{1/2+\delta} \rfloor}$, such that $G/N$ has an abelian subgroup of index bounded just in terms of $\theta$ and $\delta$.*

As it is worth to notice, the essential ingredient in the proof of Theorem 2.3 is a fundamental theorem of Breuillard, Green and Tao [1, Main Theorem].

## 3. Proofs

We start with the proof of Theorem 1.1 (A).
Let $\alpha \geq 1$ be the smallest integer such that $5^\alpha > diam(G, S)^{1/8}$. We have
$$|S^{(5^\alpha)}| \leq |G| \ll_c diam(G, S)^{1/\varepsilon} \leq diam(G, S)^{1/\varepsilon} |S| =$$
$$(diam(G, S)^{1/8})^{\frac{8}{\varepsilon}} |S| < \left(5^{8/\varepsilon}\right)^\alpha |S|,$$
and then
$$(3.1) \qquad |S^{(5^\alpha)}| < C \cdot \left(5^{8/\varepsilon}\right)^\alpha \cdot |S|$$
for some $C > 0$ which depends only on $c$ (in fact, $C = 1/c$).
Now we see that there exists $n \leq diam(G, S)^{1/4}$ such that
$$(3.2) \qquad |S^{5n}| \leq C \cdot 5^{8/\varepsilon} \cdot |S^n|.$$
Indeed, suppose that this claim does not hold. If $|S|$ remains bounded, then $diam(G, S)$ grows when $|G|$ grows (in fact, a very simple argument shows that $diam(G, S) \geq \log_{|S|} |G|$). Since $5^\alpha \leq diam(G, S)^{1/4}$ when $diam(G, S)$ is sufficiently large, it follows that
$$|S^{(5^\alpha)}| \geq C \cdot 5^{8/\varepsilon} |S^{(5^{\alpha-1})}| \geq ... \geq C^\alpha \left(5^{8/\varepsilon}\right)^\alpha |S| \geq C \left(5^{8/\varepsilon}\right)^\alpha |S|,$$
which breaks (3.1).
Let us take $N \triangleleft G$ as it is provided by Theorem 2.3 applied to the inequality (3.2) (with $\theta = C \cdot 5^{8/\varepsilon}$, $\delta = 1/4$). By definition, if $\rho : G \to G/N$ denotes the standard projection, $S^{diam(G/N, \rho(S))}$ contains a full transversal for $N$ in $G$. It follows that
$$diam(G, S) \leq diam(G, S)^{3/4} + diam(G/N, \rho(S)),$$
and then
$$|G/N| \geq diam(G/N, \rho(S)) \geq diam(G, S) - diam(G, S)^{3/4} \gg_{|S|} diam(G, S) \geq c|G|^\varepsilon.$$
Finally, let $N \triangleleft H \leqslant G$ such that $|G:H|$ is bounded and $H/N$ is abelian. We have
$$|H/H'| \geq |H/N| = \frac{|G/N|}{|G:H|} \gg_{\varepsilon,c} |G/N| \gg_{c,|S|} |G|^\varepsilon.$$

For the (easier) proof of Theorem 1.1 (B), we see how the large diameter of the abelian section influences the diameter of the full group. From Lemmas 2.1 and 2.2 we have

$$diam(G, S) \geq diam(H, \overline{S}) \geq diam(H/H', \rho(\overline{S})) \geq |H/H'|^{1/|\rho(\overline{S})|} \geq (c|G|^\varepsilon)^{1/(|G:H||S|)},$$

where $\rho : H \to H/H'$ is the standard projection. If $c \geq 1$, then obviously $c^{1/(|G:H||S|)} \geq 1$. Otherwise, we have $c^{1/(|G:H||S|)} \geq c$, and so the proof is complete.

*Remark* 3.1 (Large cyclic section). Under the hypothesis $diam(G, S) \geq c|G|^\varepsilon$, we also obtain a *cyclic* section $H/N$ of size $|G|^{O_{\varepsilon,c,|S|}(1)}$, for some $N \triangleleft H$ (and $H$ of bounded index). In fact, $H/H'$ is abelian and $|S||G : H|$-generated, so it is the direct product of at most $|S||G : H|$ cyclic factors, and at least one of these has the desired size.

The function $1 \cdot \varepsilon |S|^{-1}|G : H|^{-1}$ at the exponent in (B) is the largest possible. This can be observed looking for groups with a large abelian section and relatively small diameter.
It is not hard to see that, for every $n \geq 2$, the wreath product $G := (C_2)^n \rtimes C_n$, with respect to standard generators, has diameter $O(n) = O(\log |G|)$. Since it has a huge abelian subgroup of index $n = O(\log |G|)$, this example shows that $|G : H|^{-1}$ cannot be replaced by any larger function.
Looking at $H = G$ as test subgroup, we notice that for every finite group $G$ one has $|G/G'| \gg |G|^{1/\log |G|}$. Choosing $G$ and $S$ such that $diam(G, S) \ll \log |G|$, we see that $\varepsilon$ cannot be replaced by any larger function.
It is not hard to see that every finite group $G$ has $diam(G, S) \ll \log |G|$ for some set $S$ of size $O(\log |G|)$. Thus, $|S|^{-1}$ cannot be replaced by any larger function.
Finally, a cyclic group $G$ provides $|S| = 1$, and diameter equal to $|G|$. Thus, 1 cannot be replaced by any larger constant.

Before the proof of Proposition 1.3, we give the following interesting result.

**Lemma 3.2.** *Let $G$ be a finite group, and let $H \leqslant G$ be of minimal size among the subgroups of $G$ which provide an abelian section of maximal order. Then $H$ is nilpotent of class at most 2.*

*Proof.* We first show that $H$ is nilpotent. If not, take $M \leqslant H$ a non-normal maximal subgroup of $H$. Then $MG' = H$ or $MH' = M$, but the second is impossible because $M$ is not normal in $H$. It follows that

$$|H/H'| = |MH'/H'| = |M/M \cap H'| \leq |M/M'|.$$

Since $|M| < |H|$, then one can conclude that $H$ is nilpotent, i.e. it is the product of its Sylow's subgroups $P_1, ..., P_k$. We show that each of these has nilpotency class at most 2. Since both the derived subgroup and the center factorize through a direct product, the proof shall follow. Now $H$ has the property that $H/H'$ is larger than every other abelian section of $H$. We show that $P_i$ has the same property, for every $1 \leq i \leq k$. In fact, if not, then a too large abelian section of $H$, different from $H/H' \cong P_1/(P_1)' \times ... \times P_k/(P_k)'$, would be found. Then, the proof follows from [4, Lemma]. $\square$

*Proof of Proposition 1.3.* (i) Let us denote by $k(G)$ the number of conjugacy classes of $G$. By the orbit-counting theorem we have $k(G) = \frac{1}{|G|}\sum_{g \in G} |C_G(g)|$, where $C_G(g)$ is the centralizer of $g$ in $G$. Let $H \leqslant G$. Since $C_H(h) \subseteq C_G(h)$ for every $h \in H$, then

$$k(G) \geq \frac{1}{|G|} \sum_{h \in H} |C_H(h)| = \frac{k(H)}{|G : H|}.$$

Since $k(H)$ equals the number of irreducible characters of $H$, and $|H/H'|$ is the number of 1-dimensional characters, it follows that $k(H) \geq |H/H'|$. Now pick $H \leqslant G$ the subgroup which is produced by the point (A) of Theorem 1.1. We have

$$k(G) \geq \frac{|H/H'|}{|G : H|} \gg_{\varepsilon,c,|S|} |H/H'| \gg_{\varepsilon,c,|S|} |G|^\varepsilon.$$

(ii) Let $H \leqslant G$ be of minimal size among the subgroups of $G$ which provide an abelian section of maximal order. Then $H$ is nilpotent of class at most 2 from Lemma 3.2, and it has size at least $|G|^\varepsilon$ from Theorem 1.1 (A). $\square$

We remark that the presence of a large abelian section, as well as many conjugacy classes, does not imply a large diameter in general. For example, for prime $p \geq 2$, $G_p := SL_2(\mathbb{F}_p)$ has an abelian subgroup of size roughly $|G_p|^{1/3}$, roughly $|G_p|^{1/3}$ conjugacy classes, and still $G_p$ has very small diameter.
To conclude, we remark that normal sets of bounded size always produce a large diameter.

**Lemma 3.3** (Normal sets). *Let $G$ be a finite group, and let $S \subseteq G$ be closed under conjugation. Then*

$$diam(G, S) \gg_{|S|} |G|^{\frac{1}{|S||S|!}}.$$

*Proof.* If $S$ does not generate $G$, then the thesis holds trivially, so suppose $G = \langle S \rangle$. Let $G$ act on $S$ by conjugation, and $K := \cap_{s \in S} C_G(s)$ the kernel of this action. Certainly $Z(G) \subseteq K$. Moreover, let $k \in K$. Since $k$ centralizes every element of $S$, and every element of $G$ is a product of elements of $S$, it follows that $k \in Z(G)$. Then $Z(G) = K$, so that $G/Z(G) \cong G/K$ is permutation group over $|S|$ elements, and finally $|G/Z(G)| \leq |S|!$. From Lemmas 2.1 and 2.2 we have

$$diam(G,S) \geq diam(Z(G), \overline{S}) \geq |Z(G)|^{\frac{1}{|S||G:Z(G)|}} =$$

$$\left(\frac{|G|}{|G:Z(G)|}\right)^{\frac{1}{|S||G:Z(G)|}} \geq \left(\frac{|G|}{|S|!}\right)^{\frac{1}{|S||S|!}} \gg_{|S|} |G|^{\frac{1}{|S||S|!}}. \quad \square$$

Luca Sabatini, Dipartimento di Matematica e Informatica "Ulisse Dini",
University of Firenze, Viale Morgagni 67/a, 50134 Firenze, Italy
*Email address*: `luca.sabatini@unifi.it`